\documentclass[a4paper,11pt]{article}

\usepackage[affil-it]{authblk}
\usepackage{amsmath}
\usepackage{amsfonts}
\usepackage{amssymb}
\usepackage{graphicx}
\DeclareGraphicsExtensions{.jpg,.png,.pdf,.eps,.bmp,.gif}
\usepackage[latin1]{inputenc}
\newtheorem{theo}{Theorem}

\newtheorem{prop}[theo]{Proposition}
\newtheorem{lem}[theo]{Lemma}
\newtheorem{cor}[theo]{Corollary}

\newtheorem{rem}[theo]{Remark}

\title{On a class of self-similar processes with \\ stationary increments in higher order Wiener chaoses}
\author{Benjamin Arras
\thanks{Electronic address: \texttt{benjamin.arras@ecp.fr}}}
\affil{Ecole Centrale Paris and INRIA Regularity team\\ 
Grande Voie des Vignes, 92295~Chatenay-Malabry, France}

\begin{document}
\maketitle

\begin{abstract}
We study a class of self-similar processes with stationary increments belonging to higher order Wiener chaoses which are similar
to Hermite processes. We obtain an almost sure wavelet-like expansion of these processes.
This allows us to compute the pointwise and local Hölder regularity of sample paths and to analyse their behaviour at infinity.
We also provide some results on the Hausdorff dimension of the range and graphs of multidimensional 
anisotropic self-similar processes with stationary increments defined by multiple Wiener integrals.

\end{abstract}


\section{Introduction and background}

Self-similar processes with stationary increments (SSSI processes), {\it i.e.} processes $X$ which satisfy:
$$\forall c>0\quad \{X_{ct}:t\in\mathbb{R}_+\}\overset{(d)}{=}\{c^HX_t:t\in\mathbb{R}_+\}$$
$$\forall h>0\quad \{X_{t+h}-X_h:t\in\mathbb{R}_+\}\overset{(d)}{=}\{X_t:t\in\mathbb{R}_+\}$$
for some positive $H$, have been studied for
a long time due to their importance both in theory and in practice. Such processes appear
as limits in various normalisation procedures \cite{Lamperti,ST,MR550123}.
In addition, they are the only
possible tangent processes \cite{Falc}. In applications, they occur in various fields such as
hydrology, biomedicine and image processing. The simplest SSSI processes are simply Brownian
motion and, more generally, Lévy stable motions. Apart from these cases, the best known such
process is probably fractional Brownian motion (fBm), which was introduced in \cite{Kolmogorov} 
and popularized by \cite{MR0242239}.
A construction of SSSI processes that generalizes fBm to higher order Wiener chaoses was proposed
in \cite{MR824710}. These processes read
$$\forall t\in\mathbb{R}_+\quad X_t=\int_{\mathbb{R}^d}h_t(x_1,...,x_d)dB_{x_1}...dB_{x_d}$$
where $\{B_x:x\in\mathbb{R}\}$ is a two-sided brownian motion and where $h_t$ satisfies:

\begin{enumerate}
\item $h_t\in \hat{L}^2(\mathbb{R}^d)$, where $\hat{L}^2(\mathbb{R}^d)$ denotes the space of square-integrable symmetric functions, 
\item $\forall c>0,\quad h_{ct}(cx_1,...,cx_d)=c^{H-\frac{d}{2}}h_t(x_1,...,x_d)$,
\item $\forall \rho\geq 0, \quad h_{t+\rho}(x_1,...,x_d)-h_{t}(x_1,...,x_d)=h_{\rho}(x_1-t,...,x_d-t).$
\end{enumerate}

Properties $2.$ and $3.$ ensure self-similarity and stationarity of the increments of the process. 
Among the kernels satisfying the above properties, \cite{MR824710} considered:

\begin{itemize}
\item $\int_{0}^t||\textbf{s}^*-\textbf{x}||_2^{H-\frac{d}{2}-1}ds$,
\item $\int_{0}^t\prod_{j=1}^d(s-x_j)_+^{-(\frac{1}{2}+\frac{1-H}{d})}ds$.
\end{itemize}
where $\textbf{s}^*=(s,...,s)$, $\|.\|_2$ is the euclidian norm in $\mathbb{R}^d$ and $(x)_{+}=\max(0,x)$.\\
\\
The second kernel defines a class of processes called {\it Hermite processes}, which have been and still are the subject of considerable
interest \cite{MR550122,MR550123,MR1650532,MR2105535,Tudor}. For $d=1$, one recovers fBm (the only SSSI Gaussian process).
The process obtained with $d=2$ is called the {\it  Rosenblatt} process.

When $d>1$, the processes in the general class defined in \cite{MR824710} are somewhat difficult to analyse, because they are no longer
Gaussian. In recent years, the family of Hermite processes, and specially the Rosenblatt one, has been studied in depth from the point
of views of estimation \cite{Bard}, stochastic calculus \cite{Tudor}, distributional properties \cite{Veillette}, 
wavelet-like decomposition for $d=2$ \cite{MR2105535}, and more. Wavelet-like decompositions, in particular, are useful for investigating 
the local properties and providing synthesis algorithms. 

General results on sample paths properties of ergodic self-similar processes, such as local or pointwise regularity 
or Hausdorff dimensions were obtained in \cite{Takashima}. One such result is that, if the process has stationary increments,
then for all $t$, the pointwise Hölder exponent is equal to $H$ almost surely. In the case of fBm, a lot more is known, both
in the one and multidimensional cases \cite{MR0346925,MR618272,MR1334170,ayachejaffard2007,ayache2013}. However, in the non-Gaussian case $d>1$, there is still room for improvements. For instance, precise uniform results are not available. Of course, by stationarity of increments, 
self-similarity and the finiteness of every moments of $X_1$, one has:
$$\forall p>1,\quad\mathbb{E}[|X_t-X_s|^p]=|t-s|^{pH}\mathbb{E}[(X_1)^p].$$
By Kolmogorov's lemma, $X$ has a modification whose sample paths are Hölder-continuous of all exponents 
smaller than $H$ (we will always work with such a version in the following), but obtaining the exact 
uniform pointwise Hölder exponent is more difficult.

The aim of this article is to study a class of SSSI processes obtained by considering a particular
kernel $h$ satisfying conditions $1., 2., 3.$ above. The interest of this class is that it allows 
one to obtain a wavelet-type expansion for all values of $d$ (and not only for 
$d=1$ and $d=2$ as it is the case for Hermite processes). This expansion permits in turn to deduce 
sharp local regularity results on the sample paths. 
Our class of processes is defined by the following multiple Wiener-Itô integral representation:
\begin{equation}\label{def:process}
X^{\alpha}_t=\int_{\mathbb{R}^d}\left[||\textbf{t}^*-\textbf{x}||_2^{H-\frac{d}{2}}-||\textbf{x}||_2^{H-\frac{d}{2}}\right]dB_{x_1}...dB_{x_d}
\end{equation}
where $t\in[0,1]$, $\textbf{t}^*=(t,...,t)$ and $\alpha=H-1+\frac{d}{2}$. When $d=1$, $\{X^{\alpha}_t\}$ is a fbm.\\
\\
The main results of this article are propositions $4$, $5$ and corollary $11$. Propositions $4$ and $5$ are concerned respectively with the local regularity and the asymptotic behaviour at infinity of the process $\{X^{\alpha}_t\}$. In view of these propositions, the sample paths of the process $\{X^{\alpha}_t\}$ for $d>1$ differ significantly from the ones of the fbm since the logarithmic correction is to the power of $\frac{d}{2}$. Moreover, corollary $11$ gives us the uniform pointwise Hölder exponent (see section $4$ for a definition). It extends a well-known result regarding the fbm, namely that, almost surely, for any $t\in (0,1)$, the pointwise Hölder exponent at $t$ is equal to $H$. Nevertheless, the techniques in order to obtain this result for the considered class of processes are somewhat more involved since the collection of random variables used to define the wavelet-type expansion of $\{X^{\alpha}_t\}$ is not a collection of independent standard normal random variables (for $d>1$) anymore (see the remark after proposition $1$). Lemma $9$ is the key to obtain such an uniform result. It is based on a geometric condition which lets us obtain independence and the Carbery and Wright inequality (\cite{MR1839474}). This inequality provides us an estimate for the small ball probability of a random variable from the collection previously mentionned. It appears in the study of the convergence in law for sequences of functionals on the Wiener space (see \cite{NourdinPoly}).\\
\\
The remaining of this work is organized as follows. In the next section, we define precisely our process.
In Section \ref{sec:wd}, we prove the existence of a modification of $\{X^{\alpha}_t\}$ which is expressed as an almost-surely absolutely convergent wavelet-type series. Section \ref{sec:reg} studies the pointwise and local Hölder regularity of sample paths and the behaviour at infinity of the process. Finally, we provide in Section \ref{sec:dim} general results on the Hausdorff dimension of the range and graphs of multidimensional anisotropic SSSI processes defined by multiple Wiener integrals.

\section{Definitions and notations}\label{sec:kd}
In this section, we define the class of processes, $\{X^{\alpha}_t\}$, by means of the multiple Wiener-Itô integrals and the Riesz kernel. First of all, let us define the multiple Wiener-Itô integrals. We refer the reader to chapter 1.1.2 of \cite{MR2200233}. The multiple Wiener-Itô integrals, denoted by $I_d$, is a linear continuous application from $\hat{L}^2(\mathbb{R}^d)$ to $L^2(\Omega, \mathcal{G}, \mathbb{P})$, where $\mathcal{G}$ is the sigma-field generated by the two-sided brownian motion. Moreover, it satisfies the following properties:
\begin{itemize}
\item For all $f\in L^2(\mathbb{R}^d)$, $I_d(f)=I_d(\hat{f})$,
\item For all $f\in L^2(\mathbb{R}^d)$, $\mathbb{E}[I_d(f)]=0$,
\item For all $f\in L^2(\mathbb{R}^p)$ and $g\in L^2(\mathbb{R}^q)$,
\[
\mathbb{E}[I_p(f)I_q(g)]=
\begin{cases}
0 & p\ne q\\
p!<\hat{f};\hat{g}>_{L^2(\mathbb{R}^p)} & p=q 
\end{cases}
\]
\end{itemize}
where $\hat{f}(x_1,...,x_d)=\frac{1}{d!}\sum_{\sigma\in\mathcal{S}_d}f(x_{\sigma(1)},...,x_{\sigma(d)})$ and $\mathcal{S}_d$ is the set of permutations of $\{1,...,d\}$.\\
\\
Let us define the Riesz kernel:
\[
k_{\alpha}(\textbf{x})=\dfrac{1}{\gamma_d(\alpha)}
\begin{cases}
||\textbf{x}||_2^{\alpha-d} & \alpha-d\ne 0,2,4,6 \\
||\textbf{x}||_2^{\alpha-d}\ln\frac{1}{||\textbf{x}||_2} & \alpha-d= 0,2,4,6
\end{cases}
\]
with:
\[
\gamma_d(\alpha)=
\begin{cases}
2^{\alpha}\pi^{\frac{d}{2}}\dfrac{\Gamma(\frac{\alpha}{2})}{\Gamma(\frac{d-\alpha}{2})} & \alpha\ne d+2k,\quad \alpha\ne -2k\\
1 & \alpha=-2k\\
(-1)^{\frac{d-2}{2}}\pi^{\frac{d}{2}}2^{\alpha-1}(\dfrac{\alpha-d}{2})!\Gamma(\frac{\alpha}{2}) & \alpha=d+2k\\
\end{cases}
\]

where $k \in \mathbb{N}$ and $\Gamma(\alpha)=\int_0^{\infty}x^{\alpha-1}e^{-x}dx$.
This kernel occurs in the definition of the Riesz potential which will be used below. Let $H\in (\frac{1}{2},1)$ and set $\alpha=H+\frac{d}{2}-1$ which is, thus, in $\mathbb{R}\setminus\mathbb{N}$. Consider the following kernel:
$$h_t^{\alpha}(x_1,...,x_d)=k_{\alpha+1}(\textbf{t}^*-\textbf{x})-k_{\alpha+1}(\textbf{x})$$
where $\textbf{t}^*=(t,...,t)$. The kernel $h_t^{\alpha}$ is clearly symmetric and it is classical to show that it is a square-integrable function. Moreover, simple calculations show that Properties 2. and 3. are also verified. Thus the process:
$$\forall t\in[0,1]\quad X^{\alpha}_t=I_d(h_t^{\alpha}),$$
is an SSSI process by \cite{MR824710}.

Self-similarity and increments stationarity entail that
$$\forall (s,t)\in [0,1]^2,\quad \mathbb{E}[X^{\alpha}_tX^{\alpha}_s]=\dfrac{1}{2}\mathbb{E}[(X^{\alpha}_1)^2]\left[|t|^{2H}+|s|^{2H}-|t-s|^{2H}\right],$$
{\it i.e.} $X^{\alpha}_t$ has same covariance structure as fBm. In particular, it displays long range dependence. 
Note that, as mentioned above, Kolmogorov's lemma implies that the process has a modification whose sample paths are 
almost surely continuous. We will always work with this modification in the sequel. Using a classical patching argument, 
$X^{\alpha}$ can then be defined on the whole of $\mathbb{R}_+$.

\section{Wavelet decomposition}\label{sec:wd}

Wavelet decompositions of stochastic processes are a powerful tool for the study of sample paths properties such as regularity.  
They have been used in particular in \cite{MR2169474} and \cite{MR2508569} respectively for the cases of fractional 
Brownian sheet and linear fractional stable sheet. Following the method in these articles, we shall obtain an almost sure 
wavelet-type expansion of $\{X^{\alpha}_t\}$.  

In that view, we use a specific wavelet basis of $L^2(\mathbb{R}^d)$ (the reader is referred to \cite{MR1228209} 
for an introduction to wavelet theory). For $d=1$, we consider the Lemarié-Meyer wavelet basis,
which satisfies the following property:
\begin{itemize}
\item the scaling function $\phi$ and the mother wavelet $\psi$ are in the Schwartz space, $S(\mathbb{R})$.
\item $supp(\mathcal{F}(\phi))\subset \{\xi\in\mathbb{R}:|\xi|\leq \frac{4\pi}{3}\}$ and $supp(\mathcal{F}(\psi))\subset \{\xi\in\mathbb{R}:\frac{2\pi}{3}\leq|\xi|\leq \frac{8\pi}{3}\}$.
\item $\{2^{\frac{j}{2}}\psi(2^jx-k):(j,k)\in\mathbb{Z}\times\mathbb{Z}\}$ is an orthonormal basis of $L^2(\mathbb{R})$.
\item $\forall \beta\in\mathbb{N}\quad \int_{\mathbb{R}}x^{\beta}\psi(x)dx=0$.
\end{itemize}
The last property implies that the mother wavelet is in the Lizorkin Space, $S_0(\mathbb{R})$. 
We use the canonical construction of a multidimensional multiresolution analysis from a monodimensional 
one using tensor products. Let $E$ denote the set $\{0,1\}^d\setminus (0,...,0)$. 
For all $\epsilon\in E$, we note $\psi^{(\epsilon)}$ the function defined by:
$$\psi^{(\epsilon)}(\textbf{x})=\psi^{(\epsilon_1)}(x_1)...\psi^{(\epsilon_d)}(x_d)$$
where $\psi^{(0)}=\phi$ and $\psi^{(1)}=\psi$. Then $\{\psi^{(\epsilon)}_{j,\textbf{k}}(\textbf{x})=2^{\frac{jd}{2}}\psi^{(\epsilon)}
(2^j\textbf{x}-\textbf{k}):(j,\textbf{k},\epsilon)\in\mathbb{Z}\times\mathbb{Z}^d\times E\}$ is an orthonormal 
wavelet basis of $L^2(\mathbb{R}^d)$. By construction, for all $\epsilon\in E$, $\psi^{(\epsilon)}$ is in the 
following Lizorkin type space:
$$S_0(\mathbb{R}^d)=\{\psi\in S(\mathbb{R}^d):\forall \beta\in\mathbb{N}^d\quad \int_{\mathbb{R}^d}\textbf{x}^{\beta}\psi(\textbf{x})d\textbf{x}=0\}.$$
This space is invariant by the Riesz potential which is defined by (see \cite{MR1347689}, page $483$):
$$I^{\alpha}\psi(\textbf{x})=\int_{\mathbb{R}^d}k_{\alpha}(\textbf{x}-\textbf{y})\psi(\textbf{y})d\textbf{y}.$$
For $\psi$ in $S_0(\mathbb{R}^d)$, the following inversion formula holds:
$$I^{\alpha}\psi(\textbf{x})=\dfrac{1}{(2\pi)^d}\int_{\mathbb{R}^d}e^{i<\textbf{x},\xi>}\dfrac{\mathcal{F}(\psi)(\xi)}
{||\xi||_2^{\alpha}}d\xi.$$
We are now ready to give our wavelet-type expansion for the process $\{X^{\alpha}_t\}$. For this purpose, we proceed as follows:
\begin{itemize}
\item We show the existence of a modification of $\{X^{\alpha}_t\}$, denoted by $\{\tilde{X}^{\alpha}_t\}$, which can be expressed as a random wavelet-type series which almost surely converges absolutely at every point $t\in [0,1]$.
\item We show that this modification is continuous on $[0,1]$.
\item We conclude that these two processes are indistinguishable.
\end{itemize}

\begin{prop}\label{wavdecomp}
There exists a modification of the process $\{X^{\alpha}_t;t\in [0,1]\}$, denoted by $\{\tilde{X}^{\alpha}_t\}$, such that:
$$\forall \omega\in\Omega^*\quad\forall t\in [0,1]\quad \tilde{X}^{\alpha}_t(\omega)=\sum_{j\in\mathbb{Z}}\sum_{\textbf{k}\in\mathbb{Z}^d}\sum_{\epsilon\in 
E}2^{-jH}\left[I^{\alpha+1}(\psi^{(\epsilon)})(2^j\textbf{t}^*-\textbf{k})-I^{\alpha+1}(\psi^{(\epsilon)})
(-\textbf{k})\right]I_d(\psi^{(\epsilon)}_{j,\textbf{k}})(\omega).$$
\end{prop}
\textbf{Remark}: In \cite{MR2169474}, the wavelet expansion of the fractional brownian sheet is proved by means of the Itô-Nisio theorem (see theorem 2.1.1 in \cite{kwapien1992random}). In our case, we can not use such theorem since the collection $\{I_d(\psi^{(\epsilon)}_{j,\textbf{k}});(j,\textbf{k},\epsilon)\in\mathbb{Z}\times\mathbb{Z}^d\times E\}$ is not a collection of independent random variables. Indeed, for any $(j,\textbf{k},\textbf{k}',\epsilon,\epsilon')$, one has:
\begin{align*}
\mathbb{E}[I_d(\psi^{(\epsilon)}_{j,\textbf{k}})I_d(\psi^{(\epsilon')}_{j,\textbf{k}'})]&=d!<\hat{\psi}^{(\epsilon)}_{j,\textbf{k}};\hat{\psi}^{(\epsilon')}_{j,\textbf{k}'}>\\
&=d!<\hat{\psi}^{(\epsilon)}_{j,\textbf{k}};\psi^{(\epsilon')}_{j,\textbf{k}'}>\\
&=\sum_{\sigma\in\mathcal{S}_d}<\psi^{(\epsilon)}_{j,\textbf{k}}\circ\sigma;\psi^{(\epsilon')}_{j,\textbf{k}'}>\\
&=\sum_{\sigma\in\mathcal{S}_d}<\psi^{(\epsilon_{\sigma})}_{j,\textbf{k}_{\sigma}};\psi^{(\epsilon')}_{j,\textbf{k}'}>
\end{align*}
where $\psi\circ\sigma(x_1,...,x_d)=\psi(x_{\sigma(1)},...,x_{\sigma(d)})$, $\epsilon_{\sigma}=(\epsilon_{\sigma(1)},...,\epsilon_{\sigma(d)})$ and similarly for $\textbf{k}_{\sigma}$. Thus, if there exists $\sigma\in\mathcal{S}_d$ for which $\epsilon_{\sigma}=\epsilon'$ and $\textbf{k}_{\sigma}=\textbf{k}'$, we have $\mathbb{E}[I_d(\psi^{(\epsilon)}_{j,\textbf{k}})I_d(\psi^{(\epsilon')}_{j,\textbf{k}'})]\ne0$. For $d=3$, one can choose $\epsilon=(1,0,1)$, $\epsilon'=(1,1,0)$, $\textbf{k}=(1,2,3)$ and $\textbf{k}'=(1,3,2)$ and $\sigma=\{1,2,3\}^{\{1,3,2\}}$.\\
\\
In order to show this proposition, we need the following lemma:
\begin{lem}\label{majId}
Let $\{f_n\}$ be an orthonormal basis of $L^2(\mathbb{R}^d)$. There exists an event $\Omega^*$ of probability 
$1$ and a strictly positive random variable, $\tilde{C}_d(\omega)$, which has all its moments finite such that:
$$\forall\omega\in\Omega^*\quad\forall n\in\mathbb{N}\quad |I_d(f_n)|\leq \tilde{C}_d(\omega)(\log(e+n))^{\frac{d}{2}}.$$
\end{lem}
\textbf{Proof}: The proof follows the lines of the one of Lemma 1 in \cite{MR2027888}.
Let $a_n=(\log(e+n))^{\frac{d}{2}}$. Theorem 6.7 of \cite{MR1474726} entails that: 

$$\exists c_d>0:\quad\forall b>2, \ \forall n \in \mathbb{N},\quad\mathbb{P}\left(|I_d(f_n)|\geq ba_n\sqrt{d!}||\hat{f}_n||_{L^2(\mathbb{R}^d)}\right)\leq e^{-c_d(ba_n)^{\frac{2}{d}}} = \dfrac{1}{(e+n)^{c_db^{\frac{2}{d}}}}.$$

For $c_db^{\frac{2}{d}}>1$, $\sum \frac{1}{(e+n)^{c_db^{\frac{2}{d}}}}<\infty$, and
the Borel-Cantelli lemma implies that:
$$\exists \Omega^*,\ \mathbb{P}(\Omega^*)=1: \ \forall \omega\in\Omega^*, \ \exists n(\omega)\in\mathbb{N}: 
\ \forall n\geq n(\omega),$$
$$|I_d(f_n)(\omega)|\leq b(\log(e+n))^{\frac{d}{2}}\sqrt{d!}||\hat{f}_n||_{L^2(\mathbb{R}^d)}\leq b(\log(e+n))^{\frac{d}{2}}\sqrt{d!}.$$
Set:
\[\tau(\omega)=
\begin{cases}
\min\{n\in\mathbb{N}:k\geq n\quad |I_d(f_k)|\leq b(\log(e+k))^{\frac{d}{2}}\sqrt{d!}\}&\quad\omega\in\Omega^*\\
+\infty&\quad\omega\notin\Omega^*,
\end{cases}
\]
and define:
\[C(\omega)=
\begin{cases}
\underset{0\leq k\leq \tau(\omega)}{\sup}(|I_d(f_k)|)&\quad\omega\in\Omega^*\\
+\infty&\quad\omega\notin\Omega^*.
\end{cases}
\]
We compute the moments of any order $p>0$ of this random variable. Following the proof of Lemma 1 of \cite{MR2027888}, one gets:
\begin{eqnarray*}
\mathbb{E}(|C|^p) & = & \sum_{n=0}^{+\infty}\mathbb{E}(\underset{0\leq k\leq n}{\sup}(|I_d(f_k)|^p)\chi_{\{\tau(\omega)=n\}})\\
& \leq & \mathbb{E}(|I_d(f_0)|^p)+\sum_{n=1}^{+\infty}\sum_{k=0}^{n}\mathbb{E}(|I_d(f_k)|^p\chi_{\{|I_d(f_{n-1})|>b(\log(e+n-1))^{\frac{d}{2}}\sqrt{d!}\}})\\
& \leq & \mathbb{E}(|I_d(f_0)|^p)+\sum_{n=1}^{+\infty}\sum_{k=0}^{n}\mathbb{E}(|I_d(f_k)|^{2p})^{\frac{1}{2}}\mathbb{P}(|I_d(f_{n-1})|>b(\log(e+n-1))^{\frac{d}{2}}\sqrt{d!})^{\frac{1}{2}}.
\end{eqnarray*}

Using Theorems 5.10 and 6.7 of \cite{MR1474726}, one deduces
$$\mathbb{E}(|C|^p)\leq \mathbb{E}(|I_d(f_0)|^p)+c_p\sum_{n=1}^{+\infty}\dfrac{n+1}{(e+n-1)^{\frac{c_db^{\frac{2}{d}}}{2}}},$$
{\it i.e.}, for $\dfrac{c_db^{\frac{2}{d}}}{2}>2$,
$$\forall p>0\quad\mathbb{E}(|C|^p)<+\infty.$$
From this, one deduces that there exists a strictly positive random variable $\tilde{C}_{d}(\omega)$, 
whose every moments are finite and such that:
$$a.s.\quad\forall n\in\mathbb{N}\quad |I_d(f_n)|\leq \tilde{C}_{d}(\omega)(\log(e+n))^{\frac{d}{2}}\quad\Box$$
An indexing argument shows that the following inequality holds:
\begin{align*}
a.s.\quad\forall (j,\textbf{k},\epsilon)\in\mathbb{Z}\times\mathbb{Z}^d\times E\quad |I_d(\psi^{(\epsilon)}_{j,\textbf{k}})|&\leq \tilde{C}_{1,d}(\omega)(\log(e+|j|+||\textbf{k}||_1))^{\frac{d}{2}}\\
&\leq \tilde{C}_{1,d}(\omega)(\log(e+|j|))^{\frac{d}{2}}(\log(e+||\textbf{k}||_1))^{\frac{d}{2}},
\end{align*}
where we have used the two trivial facts:
\begin{itemize}
\item the function $x\rightarrow x^{\frac{d}{2}}$ is increasing on $\mathbb{R}_+$
\item $\forall (x,y)\in \mathbb{R}_+\quad \log(e+x+y)\leq \log(e+x)\log(e+y)$
\end{itemize}

Let us now move to the proof of Proposition \ref{wavdecomp}. It is similar to that of the wavelet-type 
expansion of fBm in \cite{MR1755100} and in section 4.2 of \cite{ayache2013}.\\
\\
\textbf{Proof}: Let $t$ be in $[0,1]$. Since $h_t^{\alpha}$ is in $\hat{L}^2(\mathbb{R}^d)$, we can expand it into the wavelet basis $\{\psi^{(\epsilon)}_{j,\textbf{k}}\}$. Using Itô isometry we have, in $(L^2)$:
$$X^{\alpha}_t=\sum_{j\in\mathbb{Z}}\sum_{\textbf{k}\in\mathbb{Z}^d}\sum_{\epsilon\in E}<h_t^{\alpha},\psi^{(\epsilon)}_{j,\textbf{k}}>I_d(\psi^{(\epsilon)}_{j,\textbf{k}}).$$
Let us show that the right hand side converges absolutely for every $t$ in $[0,1]$ on $\Omega^*$. It suffices to show that:
$$\forall\omega\in \Omega^*\quad \sum_{j\in\mathbb{Z}}\sum_{\textbf{k}\in\mathbb{Z}^d}\sum_{\epsilon\in E}|<h_t^{\alpha},\psi^{(\epsilon)}_{j,\textbf{k}}>||I_d(\psi^{(\epsilon)}_{j,\textbf{k}})|<+\infty$$
For this purpose, write:
$$<h_t^{\alpha},\psi^{(\epsilon)}_{j,\textbf{k}}>=\int_{\mathbb{R}^d}[k_{\alpha+1}(\textbf{t}^*-\textbf{x})-k_{\alpha+1}(-\textbf{x})]\psi^{(\epsilon)}_{j,\textbf{k}}(\textbf{x})d\textbf{x}.$$
The change of variables $2^j\textbf{x}-\textbf{k}=\textbf{u}$ yields:
$$<h_t^{\alpha},\psi^{(\epsilon)}_{j,\textbf{k}}>=2^{-jH}\int_{\mathbb{R}^d}[k_{\alpha+1}(2^j\textbf{t}^*-\textbf{k}-\textbf{u})-k_{\alpha+1}(-\textbf{k}-\textbf{u})]\psi^{(\epsilon)}(\textbf{u})d\textbf{u}.$$
Split the integral into two terms to get:
$$<h_t^{\alpha},\psi^{(\epsilon)}_{j,\textbf{k}}>=2^{-jH}\left[I^{\alpha+1}(\psi^{(\epsilon)})(2^j\textbf{t}^*-\textbf{k})-I^{\alpha+1}(\psi^{(\epsilon)})(-\textbf{k})\right].$$
We need to distinguish between $j\leq -1$ and $j\geq 0$. Since $I^{\alpha+1}(\psi^{(\epsilon)})\in S_0(\mathbb{R}^d)$, the function $f$ defined on $\mathbb{R}$ by:
$$f(u)=I^{\alpha+1}(\psi^{(\epsilon)})(\textbf{u}^*-\textbf{k})$$
is in $C^{\infty}(\mathbb{R})$. Thus, by the mean value theorem:
$$\exists v\in [0,2^jt]:\quad I^{\alpha+1}(\psi^{(\epsilon)})(2^j\textbf{t}^*-\textbf{k})-I^{\alpha+1}(\psi^{(\epsilon)})(-\textbf{k})=2^jt\dfrac{df}{du}(v).$$
By the chain rule,
$$\dfrac{df}{du}(v)=\sum_{i=1}^{d}\dfrac{\partial}{\partial u_i}(I^{\alpha+1}(\psi^{(\epsilon)}))(\textbf{v}^*-\textbf{k}).$$
Using the rapid decrease of $\frac{\partial}{\partial u_i}(I^{\alpha+1}(\psi^{(\epsilon)}))$, one gets:
$$\forall p>0\quad |\dfrac{df}{du}(v)|\leq \dfrac{C_{p,d}}{(1+||\textbf{v}^*-\textbf{k}||_2)^p}$$
For $j\leq -1$, $v\in[0,2^jt]\subset [0,\frac{1}{2}]$. Thus, for $||\textbf{k}||_2$ large enough,
the following estimate holds:
$$\forall t\in [0,1]\quad |<h_t^{\alpha},\psi^{(\epsilon)}_{j,\textbf{k}}>|\leq 2^{j(1-H)}\dfrac{C_{p,d}}{(1+||\textbf{k}||_2)^p}.$$
By Lemma \ref{majId},
$$\forall\omega\in\Omega^*\quad\forall (j,\textbf{k},\epsilon)\in\mathbb{Z}\times\mathbb{Z}^d\times E\quad |I_d(\psi^{(\epsilon)}_{j,\textbf{k}})(\omega)|\leq \tilde{C}(\omega)(\log(e+|j|))^{\frac{d}{2}}(\log(e+||\textbf{k}||_1))^{\frac{d}{2}}.$$
Thus,
$$\|<h_.^{\alpha},\psi^{(\epsilon)}_{j,\textbf{k}}>\|_{\infty,[0,1]}|I_d(\psi^{(\epsilon)}_{j,\textbf{k}})(\omega)|\leq C'_{p,d}(\omega)\dfrac{2^{j(1-H)}}{(1+||\textbf{k}||_2)^p}(\log(e+|j|))^{\frac{d}{2}}(\log(e+||\textbf{k}||_1))^{\frac{d}{2}}.$$
Since $H<1$, almost sure uniform convergence for the ``low frequency'' part $j<0$ for $p$ sufficiently 
large does hold. The high frequency part $j\geq 0$ requires a bit more work. The inversion formula 
for the Riesz potential yields:
$$I^{\alpha+1}(\psi^{(\epsilon)})(2^j\textbf{t}^*-\textbf{k})=\dfrac{1}{(2\pi)^d}\int_{\mathbb{R}^d}e^{i<2^j\textbf{t}^*-\textbf{k},\xi>}\dfrac{\mathcal{F}(\psi^{(\epsilon)})(\xi)}{||\xi||_2^{\alpha+1}}d\xi.$$
Making $p$ integrations by parts in the direction $i$ where $||2^j\textbf{t}^*-\textbf{k}||_2\leq C_{d}|2^jt-k_i|$, one gets:
$$|I^{\alpha+1}(\psi^{(\epsilon)})(2^j\textbf{t}^*-\textbf{k})|\leq \dfrac{C'_{p,d}}{(1+||2^j\textbf{t}^*-\textbf{k}||_2)^p}.$$
Thus, for $j\geq 0$:
$$|<h_t^{\alpha},\psi^{(\epsilon)}_{j,\textbf{k}}>|\leq C''_{p,d}2^{-jH}\left[\dfrac{1}{(1+||2^j\textbf{t}^*-\textbf{k}||_2)^p}+\dfrac{1}{(1+||\textbf{k}||_2)^p}\right].$$
In order to conclude, we need to estimate the following quantity:
$$\sum_{\textbf{k}\in\mathbb{Z}^d}\dfrac{(\log(e+||\textbf{k}||_1))^{\frac{d}{2}}}{(1+||2^j\textbf{t}^*-\textbf{k}||_2)^p}.$$
Denote $\lfloor2^j\textbf{t}\rfloor^*=(\lfloor2^jt\rfloor,...,\lfloor2^jt\rfloor)\in\mathbb{Z}^d$ where $\lfloor u\rfloor$ is the integer part of the real $u$. 
A change of variable in the previous sum yields:
$$\sum_{\textbf{k}\in\mathbb{Z}^d}\dfrac{(\log(e+||\textbf{k}||_1))^{\frac{d}{2}}}{(1+||2^j\textbf{t}^*-\textbf{k}||_2)^p}=\sum_{\textbf{k}\in\mathbb{Z}^d}\dfrac{(\log(e+||\textbf{k}+\lfloor2^j\textbf{t}\rfloor^*||_1))^{\frac{d}{2}}}{(1+||2^j\textbf{t}^*-\lfloor2^j\textbf{t}\rfloor^*-\textbf{k}||_2)^p}.$$
By the triangular inequality, 
$$\log(e+||\textbf{k}+\lfloor2^j\textbf{t}\rfloor^*||_1)\leq \log(e+||\textbf{k}||_1+||\lfloor2^j\textbf{t}\rfloor^*||_1)\leq C' \log(e+||\textbf{k}||_1)\log(e+||\lfloor2^j\textbf{t}\rfloor^*||_1).$$
Consequently, one obtains:
$$\sum_{\textbf{k}\in\mathbb{Z}^d}\dfrac{(\log(e+||\textbf{k}||_1))^{\frac{d}{2}}}{(1+||2^j\textbf{t}^*-\textbf{k}||_2)^p}\leq C'(\log(e+||\lfloor2^j\textbf{t}\rfloor^*||_1))^{\frac{d}{2}}\underset{\textbf{z}\in[0,1]^d}{\sup}\left(\sum_{\textbf{k}\in\mathbb{Z}^d}\dfrac{(\log(e+||\textbf{k}||_1))^{\frac{d}{2}}}{(1+||\textbf{z}-\textbf{k}||_2)^p}\right).$$
Thus, for $j\geq 0$:
$$\sum_{\textbf{k}\in\mathbb{Z}^d}|<h_t^{\alpha},\psi^{(\epsilon)}_{j,\textbf{k}}>||I_d(\psi^{(\epsilon)}_{j,\textbf{k}})(\omega)|\leq C'''_{p,d}(\omega)2^{-jH}(\log(e+d2^jt))^{\frac{d}{2}}(\log(e+|j|))^{\frac{d}{2}},$$
which concludes the proof .$\Box$\\
\\
\textbf{Remark}: This wavelet expansion does not extend straightforwardly to every stochastic process represented by $I_d(h_t)$ where $h_t$ satisfies properties 1.,2. and 3. Indeed, for the hermite processes, the wavelet coefficients of the kernel $h_t$ are equal to:
\begin{align*}
<h_t;\psi^{(\epsilon)}_{j,\textbf{k}}>=&c(H_0)\int_0^t\int_{\mathbb{R}^d}\bigg(\prod_{j=1}^d\dfrac{(s-x_j)^{H_0-1}_+}{\Gamma(H_0)}\bigg)\psi^{(\epsilon)}_{j,\textbf{k}}(x_1,...,x_d)dx_1...dx_dds\\
&=c(H_0)\int_0^t(I_{+}^{H_0}\otimes...\otimes I_{+}^{H_0})(\psi^{(\epsilon)}_{j,\textbf{k}})(\textbf{s}^*)ds\\
&=c(H_0)\int_0^tI_{+}^{H_0}(\psi^{(\epsilon_1)}_{j,k_1})(s)...I_{+}^{H_0}(\psi^{(\epsilon_d)}_{j,k_d})(s)ds\\
&=c(H_0)2^{-jH}\int_0^{2^jt}I_{+}^{H_0}(\psi^{(\epsilon_1)})(s-k_1)...I_{+}^{H_0}(\psi^{(\epsilon_d)})(s-k_d)ds
\end{align*}
where $H_0=\frac{H}{d}+\frac{d-2}{2d}$, $I_{+}^{H_0}(\psi)(t)=\frac{1}{\Gamma(H_0)}\int_{-\infty}^t(t-s)^{H_0-1}_{+}\psi(s)ds$ is the fractional integral of Weyl type of order $H_0$ (see \cite{MR1347689}) and $c(H_0)$ a certain normalizing constant. Take for instance, $d=2$, $\epsilon=(1,0)$ and $(j,\textbf{k})\in\mathbb{Z}\times\mathbb{Z}^2$. We obtain:
\begin{align*}
<h_t;\psi^{(\epsilon)}_{j,\textbf{k}}>=c(H_0)2^{-jH}\int_0^{2^jt}I_{+}^{H_0}(\psi)(s-k_1)I_{+}^{H_0}(\phi)(s-k_2)ds
\end{align*}
But the scaling function $\phi$ is not in the Lizorkin space $S_0(\mathbb{R})$. Thus, $I_{+}^{H_0}(\phi)$ does not decrease sufficiently fast at infinity. Actually, by Lemma $1$ of \cite{MR2105535}, we have:
\begin{align*}
\exists C>0\quad\forall x\in\mathbb{R}\quad|I_{+}^{H_0}(\phi)(x)|\leq \dfrac{C}{(1+|x|)^{1-H_0}}.
\end{align*}
We postpone the proof of the continuity of the modification, $\{\tilde{X}^{\alpha}_t\}$, to the next section. Indeed in proposition 4, we find a modulus of continuity for $\{\tilde{X}^{\alpha}_t\}$ using its definition and consequently, $\{\tilde{X}^{\alpha}_t\}$ is continuous on $[0,1]$. Therefore, in the sequel, by indistinguishability, we use $\{X^{\alpha}_t\}$ in order to designate both processes.

\section{Uniform and Local Regularity of the sample paths}\label{sec:reg}

From Proposition \ref{wavdecomp}, one sees that the process $X^{\alpha}_t$ may be written:
$$X^{\alpha}_t=X^{l,\alpha}_t+{X}^{h,\alpha}_t,$$
where $X^{l,\alpha}_t$ is the low frequency part and $X^{h,\alpha}_t$ the high frequency part, namely:
$$X^{l,\alpha}_t=\sum_{j\in\mathbb{Z}_-}\sum_{\textbf{k}\in\mathbb{Z}^d}\sum_{\epsilon\in E}2^{-jH}\left[I^{\alpha+1}(\psi^{(\epsilon)})(2^j\textbf{t}^*-\textbf{k})-I^{\alpha+1}(\psi^{(\epsilon)})(-\textbf{k})\right]I_d(\psi^{(\epsilon)}_{j,\textbf{k}}),$$
$$X^{h,\alpha}_t=\sum_{j\in\mathbb{N}}\sum_{\textbf{k}\in\mathbb{Z}^d}\sum_{\epsilon\in E}2^{-jH}\left[I^{\alpha+1}(\psi^{(\epsilon)})(2^j\textbf{t}^*-\textbf{k})-I^{\alpha+1}(\psi^{(\epsilon)})(-\textbf{k})\right]I_d(\psi^{(\epsilon)}_{j,\textbf{k}}).$$
In order to study the regularity of the stochastic process, we will show that the low frequency part is $C^{\infty}(\mathbb{R})$ and then focus on the regularity of the high frequency part.
\begin{prop}\label{basf}
Almost surely, the process $\{X^{l,\alpha}_t\}$ is in $C^{\infty}(\mathbb{R})$.\\
\end{prop}
\textbf{Proof}: Let $g_{j,\textbf{k}}$ be defined by:
$$\forall t\in[0,1]\quad g_{j,\textbf{k}}(t)=I^{\alpha+1}(\psi^{(\epsilon)})(2^j\textbf{t}^*-\textbf{k}).$$
We need to show that the following random series converges uniformly on $[0,1]$ almost surely:
$$\sum_{j\in\mathbb{Z}_-}\sum_{\textbf{k}\in\mathbb{Z}^d}\sum_{\epsilon\in E}2^{-jH}\dfrac{d^mg_{j,\textbf{k}}}{dt^m}(t)I_d(\psi^{(\epsilon)}_{j,\textbf{k}}), \quad m\in \mathbb{N}^*.$$
The chain rule formula entails that:
$$\dfrac{d^mg_{j,\textbf{k}}}{dt^m}(t)=2^{jm}\sum_{(i_1,...,i_m)\in\{1,...,d\}^m}\dfrac{\partial^m}{\partial x_{i_1}...\partial x_{i_m}}(I^{\alpha+1}(\psi^{(\epsilon)}))(2^j\textbf{t}^*-\textbf{k}).$$
Since $I^{\alpha+1}(\psi^{(\epsilon)})\in S_0(\mathbb{R}^d)$, one obtains the following estimate:
$$\forall p>0\quad |\dfrac{d^mg_{j,\textbf{k}}}{dt^m}(t)|\leq 2^{jm}\dfrac{C_{p,d,m}}{(1+||2^j\textbf{t}^*-\textbf{k}||_2)^p}.$$
Moreover, since $t\in[0,1]$ and $j\leq -1$, $2^j\textbf{t}^*\in[0,\frac{1}{2}]^d$. 
Thus, for $||\textbf{k}||_2$ sufficiently large, one has:
$$|\dfrac{d^mg_{j,\textbf{k}}}{dt^m}(t)|\leq 2^{jm}\dfrac{C_{p,d,m}}{(1+||\textbf{k}||_2)^p}.$$
Finally, since $H<1\leq m$ and using Lemma \ref{majId}, one obtains the almost sure uniform convergence of the random series. $\Box$\\
\\
The next two results are analogues of Theorems 1 and 2 of \cite{MR2169474}. The first one gives a modulus of
continuity for the process $X^{\alpha}_t$ while the second one deals with the asymptotic behaviour of the 
process as $t\rightarrow+\infty$. Their proofs are also similar.
\begin{prop}\label{locholder}
There exists a strictly positive random variable $A_d$ of finite moments of any order and a constant, $b_d>1$, such that:
$$\forall\omega\in\Omega^*\quad \underset{(s,t)\in[0,1]}{\sup}\dfrac{|X^{\alpha}_t(\omega)-X^{\alpha}_s(\omega)|}{|t-s|^H(\log(b_d+|t-s|^{-1}))^{\frac{d}{2}}}\leq A_d(\omega).$$
\end{prop}
This result should be compared with Remark 2.3 p.166 and the sentence that follows its proof on p.167 
of \cite{Takashima}, which provide a related (but weaker) result, namely that the local Hölder 
exponent of $X^{\alpha}$ (as well as the ones of all
Hermite processes) is equal to $H$: the above proposition entails that the local Hölder exponent 
of $X^{\alpha}$ is not smaller than $H$, and Theorem \ref{expchamp} below will ensure that 
it is not greater than $H$. See also Lemma 5.9 of \cite{MR2759163}, which provides a similar result for Hermite processes.\\
\\
\textbf{Proof}: Since $t\rightarrow X^{l,\alpha}_t$ is in $C^{\infty}(\mathbb{R})$ almost surely, we need only to deal with
the high frequency component $X^{h,\alpha}_t$. Let $\omega\in\Omega^*$ and $(s,t)$ be in $[0,1]$. 
There exists a positive integer $j_0$ such that:
$$2^{-j_0-1}<|t-s|\leq 2^{-j_0}.$$
Write:
$$X^{h,\alpha}_t-X^{h,\alpha}_s=(I)+(II),$$
where
$$(I)=\sum_{j=0}^{j_0}\sum_{\textbf{k}\in\mathbb{Z}^d}\sum_{\epsilon\in E}2^{-jH}\left[I^{\alpha+1}(\psi^{(\epsilon)})(2^j\textbf{t}^*-\textbf{k})-I^{\alpha+1}(\psi^{(\epsilon)})(2^j\textbf{s}^*-\textbf{k})\right]I_d(\psi^{(\epsilon)}_{j,\textbf{k}})$$
and
$$(II)=\sum_{j\geq j_0+1}\sum_{\textbf{k}\in\mathbb{Z}^d}\sum_{\epsilon\in E}2^{-jH}\left[I^{\alpha+1}(\psi^{(\epsilon)})(2^j\textbf{t}^*-\textbf{k})-I^{\alpha+1}(\psi^{(\epsilon)})(2^j\textbf{s}^*-\textbf{k})\right]I_d(\psi^{(\epsilon)}_{j,\textbf{k}}).$$
Let us consider the first sum. Using Lemma \ref{majId} and following the lines of the proof of 
Theorem 1 of \cite{MR2169474}, there exists a deterministic constant $b_{1,d}>1$ such that:
$$|(I)|\leq A_{1,d}(\omega)|t-s|\sum_{j=0}^{j_0}2^{(1-H)j}(b_{1,d}+j)^{\frac{d}{2}}\leq A'_{1,d}(\omega)|t-s|2^{(j_0+1)(1-H)}(b_{1,d}+j_0)^{\frac{d}{2}}.$$
Similarly, for the second sum, there exists a deterministic constant $b_{2,d}>1$ such that:
$$|(II)|\leq A_{2,d}(\omega)\sum_{j=j_0+1}^{\infty}2^{-jH}(b_{2,d}+j)^{\frac{d}{2}}\leq A'_{2,d}(\omega)2^{-j_0H}(b_{2,d}+j_0)^{\frac{d}{2}}.$$
Putting the results altogether and using the definition of $j_0$ yields the desired conclusion. $\Box$\\

\begin{prop}\label{infi} 
There exists a strictly positive random variable $B_d$ of finite moments of any order and a constant 
$c_d>3$, such that:
$$\forall\omega\in\Omega^*\quad \underset{t\in\mathbb{R_+}}{\sup}\dfrac{|X^{\alpha}_t(\omega)|}{(1+|t|)^H(\log\log(c_d+|t|))^{\frac{d}{2}}}\leq B_d(\omega).$$
\end{prop}
\textbf{Proof}: The proof is quite similar to the previous one except that the behaviour at infinity 
of the process will be governed by the low frequency part.
Let $t$ be a real strictly greater than $1$ and let $j_1$ denote the integer such that:
$$2^{j_1}\leq t<2^{j_1+1}.$$
Let $\omega\in \Omega^*$. We consider the following quantity that we split into three parts:
$$\sum_{j\in\mathbb{Z}}\sum_{\textbf{k}\in\mathbb{Z}^d}\sum_{\epsilon\in E}2^{-jH}\left[I^{\alpha+1}(\psi^{\epsilon})(2^j\textbf{t}^*-\textbf{k})-I^{\alpha+1}(\psi^{\epsilon})(-\textbf{k})\right](\log(e+|j|+||\textbf{k}||))^{\frac{d}{2}}=(I)+(II)+(III)$$
where,
$$(I)=\sum_{j=-\infty}^{-(j_1+1)}\sum_{\textbf{k}\in\mathbb{Z}^d}\sum_{\epsilon\in E}2^{-jH}\left[I^{\alpha+1}(\psi^{\epsilon})(2^j\textbf{t}^*-\textbf{k})-I^{\alpha+1}(\psi^{\epsilon})(-\textbf{k})\right](\log(e+|j|+||\textbf{k}||))^{\frac{d}{2}},$$
$$(II)=\sum_{j=-j_1}^{-1}\sum_{\textbf{k}\in\mathbb{Z}^d}\sum_{\epsilon\in E}2^{-jH}\left[I^{\alpha+1}(\psi^{\epsilon})(2^j\textbf{t}^*-\textbf{k})-I^{\alpha+1}(\psi^{\epsilon})(-\textbf{k})\right](\log(e+|j|+||\textbf{k}||))^{\frac{d}{2}},$$
$$(III)=\sum_{j=0}^{\infty}\sum_{\textbf{k}\in\mathbb{Z}^d}\sum_{\epsilon\in E}2^{-jH}\left[I^{\alpha+1}(\psi^{\epsilon})(2^j\textbf{t}^*-\textbf{k})-I^{\alpha+1}(\psi^{\epsilon})(-\textbf{k})\right](\log(e+|j|+||\textbf{k}||))^{\frac{d}{2}}.$$
Let us deal first with the high frequency part. We have:
$$|(III)|\leq h_3(t)+h_3(0)$$
where $h_3$ is the function defined by:
$$h_3(t)=\sum_{j=0}^{\infty}\sum_{\textbf{k}\in\mathbb{Z}^d}\sum_{\epsilon\in E}2^{-jH}|I^{\alpha+1}(\psi^{\epsilon})(2^j\textbf{t}^*-\textbf{k})|(\log(e+|j|+||\textbf{k}||))^{\frac{d}{2}}.$$
Reasoning as above, one gets:
$$h_3(t)\leq C_3\sum_{j=0}^{+\infty}2^{-jH}\left(\log(e+j+||2^j\textbf{t}^*||)\right)^{\frac{d}{2}}.$$
Using that:
$$\forall a,b>0\quad \log(e+a+b)\leq \log(e+a)+\log(e+b),$$
and, for $d>2$, by a convexity argument:
$$\forall a,b>0\quad (a+b)^{\frac{d}{2}}\leq c_d(a^{\frac{d}{2}}+b^{\frac{d}{2}}),$$
one gets:
$$h_3(t) \leq C_{3d}\sum_{j=0}^{+\infty}2^{-jH}\left(\log(e+j)\right)^{\frac{d}{2}}+C_{3d}\sum_{j=0}^{+\infty}2^{-jH}\left(\log(e+2^j||\textbf{t}^*||)\right)^{\frac{d}{2}}.$$
Finally, since $\log(e+2^j||\textbf{t}^*||)\leq \log(e+2^j)+\log(e+||\textbf{t}^*||)$:
$$h_3(t)\leq C'_{3d} \left(\log(e+||\textbf{t}^*||)\right)^{\frac{d}{2}}.$$
Consider now the second term. One has:
$$|(II)|\leq h_2(t)+h_2(0),$$
where $h_2$ is the function defined by:
$$h_2(t)=\sum_{j=-j_1}^{-1}\sum_{k\in\mathbb{Z}^d}\sum_{\epsilon\in E}2^{-jH}|I^{\alpha+1}(\psi^{(\epsilon)})(2^j\textbf{t}^*-\textbf{k})|\left(\log(e+|j|+||\textbf{k}||)\right)^{\frac{d}{2}}.$$
It is clear that:
$$h_2(t)\leq C_2\sum_{j=1}^{j_1}2^{jH}\left(\log(e+j+2^{-j}||\textbf{t}^*||)\right)^{\frac{d}{2}}.$$
Using the definition of $j_1$ and elementary computations, one obtains:
$$h_2(t)\leq C_{2,d} (1+|t|)^H\left(\log\log(c_d+|t|)\right)^{\frac{d}{2}}.$$
Similarly for $h_2(0)$, one has:
$$h_2(0)\leq C'_{2,d} (1+|t|)^H\left(\log\log(c_d+|t|)\right)^{\frac{d}{2}}$$
For the last term, we apply the mean value theorem to the function $f$ defined above to get:
$$|(I)|\leq C_{1,d}\sum_{j=-\infty}^{-(j_1+1)}\sum_{\textbf{k}\in\mathbb{Z}^d}\sum_{\epsilon\in E}2^{(1-H)j}|t|\dfrac{\left(\log(e+|j|+||\textbf{k}||)\right)^{\frac{d}{2}}}{(1+||\textbf{k}||)^p}\leq C'_{1,d}|t|,$$
and the result follows. $\Box$\\

\begin{rem}
In order to study the local regularity of the process, we will only need the behaviour at infinity of the high frequency part.
\end{rem}
We now compute the uniform almost sure pointwise Hölder exponent of $X^{\alpha}$. Since the paths of the process
are nowhere differentiable, it is by definition equal to:
$$\gamma_{X^{\alpha}}(t)=\sup\{\gamma>0:\quad\underset{\rho\rightarrow 0}{\limsup}\dfrac{|X^{\alpha}_{t+\rho}-X^{\alpha}_{t}|}{|\rho|^{\gamma}}<+\infty\}.$$
Proposition \ref{locholder} entails that, almost surely:
$$\forall t\in(0,1)\quad\gamma_{X^{\alpha}}(t)\geq H.$$
We will show the converse inequality. In that view, we introduce a stochastic field 
$\textbf{X}^{\alpha}_{\textbf{t}}$ that is closely
related to $X^{\alpha}_t$ and defined as follows:
$$\forall\omega\in\Omega^*, \ \forall \textbf{t}\in\mathbb{R}^d\quad\textbf{X}^{\alpha}_{\textbf{t}}(\omega)=\sum_{j\in\mathbb{N}}\sum_{\textbf{k}\in\mathbb{Z}^d}\sum_{\epsilon\in E}2^{-jH}\left[I^{\alpha+1}(\psi^{(\epsilon)})(2^j\textbf{t}-\textbf{k})-I^{\alpha+1}(\psi^{(\epsilon)})(-\textbf{k})\right]I_d(\psi^{\epsilon}_{j,\textbf{k}}).$$
It is easily checked that $\textbf{X}^{\alpha}_{\textbf{t}}$ is well defined and that its sample paths are 
continuous. Moreover, proceeding as in the proof of Proposition 5 for the high frequency part (III), one can prove that:
$$a.s., \ \forall \textbf{t}\in \mathbb{R}^d, \ |\textbf{X}^{\alpha}_{\textbf{t}}|\leq \tilde{C}_d(\omega)\left(\log(e+||\textbf{t}||)\right)^{\frac{d}{2}}.$$
Define the following operator on $S_0(\mathbb{R}^d)$:
$$\forall \textbf{x}\in\mathbb{R}^d\quad D^{\alpha+1}(\psi)(\textbf{x})=\frac{1}{(2\pi)^d}\int_{\mathbb{R}^d}e^{i<\textbf{x};\xi>}||\xi||^{\alpha+1}\mathcal{F}(\psi)(\xi)d\xi.$$
The Fourier transform of the dilated and translated version of this function is equal to:
$$\forall\xi\in\mathbb{R}^d\quad \mathcal{F}(D^{\alpha+1}(\psi)_{j,\textbf{k}})(\xi)=2^{-jd}e^{-i2^{-j}<\textbf{k};\xi>}\mathcal{F}(D^{\alpha+1}(\psi))(2^{-j}\xi),$$
where $\mathcal{F}(D^{\alpha+1}(\psi))(\xi)=||\xi||^{\alpha+1}\mathcal{F}(\psi)(\xi)$.

\begin{lem}\label{Idchamp} 
Almost surely,
$$\forall (j,\textbf{k},\epsilon)\mathbb{N}\times\mathbb{Z}^d\times E\quad I_{d}(\psi^{(\epsilon)}_{j,\textbf{k}})=2^{j(H+d)}\int_{\mathbb{R}^d}\textbf{X}^{\alpha}_{\textbf{t}}D^{\alpha+1}(\psi^{\epsilon})(2^j\textbf{t}-\textbf{k})d\textbf{t}.$$
\end{lem}
\textbf{Proof}: For $\omega\in \Omega^*$ and $(j,\textbf{k},\epsilon)\in\mathbb{N}\times\mathbb{Z}^d\times E$, we compute the following quantity:
$$A:=\int_{\mathbb{R}^d}\textbf{X}^{\alpha}_{t_1,...,t_d}D^{\alpha+1}(\psi^{(\epsilon)})(2^jt_1-k_1,...,2^jt_d-k_d)dt_1...dt_d.$$
Since the function $D^{\alpha+1}(\psi^{(\epsilon)})$ is in $S_0(\mathbb{R}^d)$, the asymptotic growth 
properties of the stochastic field $\{\textbf{X}^{\alpha}_t\}$ allows one to invert the sum and the integral, so that:
\begin{eqnarray*}
A & = & \sum_{j',\textbf{k}',\epsilon'}2^{-j'H}\int_{\mathbb{R}^d}I^{\alpha+1}(\psi^{(\epsilon')})(2^{j'}\textbf{t}-\textbf{k}')D^{\alpha+1}(\psi^{(\epsilon)})(2^j\textbf{t}-\textbf{k})d\textbf{t} \
I_d(\psi^{(\epsilon')}_{j',\textbf{k}'})\\
& =: & \sum_{j',\textbf{k}',\epsilon'} A_{j',\textbf{k}',\epsilon'} \ I_d(\psi^{(\epsilon')}_{j',\textbf{k}'}).
\end{eqnarray*}
Let us compute the inner products $A_{j',\textbf{k}',\epsilon'}$. By Fourier isometry,
$$A_{j',\textbf{k}',\epsilon'}=2^{-j'H}\dfrac{1}{(2\pi)^d}\int_{\mathbb{R}^d}2^{-j'd}e^{-i2^{-j'}<\textbf{k}';\xi>}\mathcal{F}(I^{\alpha+1}(\psi^{(\epsilon')}))(2^{-j'}\xi)2^{-jd}e^{+i2^{-j}<\textbf{k};\xi>}\overline{\mathcal{F}(D^{\alpha+1}(\psi^{(\epsilon)}))(2^{-j}\xi)} \ d\xi.$$
Using the inversion formula and the definition of the operator $D^{\alpha+1}$, one gets:
$$A_{j',\textbf{k}',\epsilon'}=\dfrac{2^{-j(H+d)}}{(2\pi)^d}\int_{\mathbb{R}^d}2^{-j\frac{d}{2}}2^{-j'\frac{d}{2}}e^{-i2^{-j'}<\textbf{k}';\xi>}e^{+i2^{-j}<\textbf{k};\xi>}\mathcal{F}(\psi^{(\epsilon')})(2^{-j'}\xi)\overline{\mathcal{F}(\psi^{(\epsilon)})(2^{-j}\xi)}d\xi.$$
Since:
$$\mathcal{F}(\psi^{(\epsilon)}_{j,\textbf{k}})(\xi)=2^{-j\frac{d}{2}}e^{-i2^{-j}<\textbf{k};\xi>}\mathcal{F}(\psi^{(\epsilon)})(2^{-j}\xi),$$
one finally obtains,
$$A_{j',\textbf{k}',\epsilon'}=2^{-j(H+d)}\int_{\mathbb{R}^d}\psi^{(\epsilon)}_{j,\textbf{k}}(\textbf{x})\psi^{(\epsilon')}_{j',\textbf{k}'}(\textbf{x})dx.$$
But $\{\psi^{(\epsilon)}_{j,\textbf{k}}\}$ is an orthonormal basis of $L^2(\mathbb{R}^d)$, so that:
$$\int_{\mathbb{R}^d}\textbf{X}^{\alpha}_{\textbf{t}}D^{\alpha+1}(\psi^{\epsilon})(2^j\textbf{t}-\textbf{k})d\textbf{t}=2^{-j(H+d)}I_{d}(\psi^{(\epsilon)}_{j,\textbf{k}}).\quad \Box$$
\\
We will show that the pointwise Hölder exponent of the stochastic field $\textbf{X}^{\alpha}_{\textbf{t}}$ is 
everywhere equal to $H$ almost surely. As a first step, we will exhibit an event of probability one such that,
for any $\textbf{t}_0\in[0,1]^d$, there exists a random variable $I_{d}(\psi^{(\epsilon)}_{j,\textbf{k}})$ 
which is both ``large'' enough and whose index is ``close'' to $\textbf{t}_0$. The construction of such an event is similar to lemma 4.1 in \cite{ayachejaffard2007} but, as mentionned in the introduction, somewhat more involved due to the fact that the collection $\{I_d(\psi_{j,\textbf{k}}^{(\epsilon)});(j,\textbf{k},\epsilon)\in\mathbb{Z}\times\mathbb{Z}^d\times E\}$ is not a collection of independent standard normal random variables for $d>1$. In that view, we will use the following result from \cite{MR1839474}:

\begin{lem}\label{Carb}\cite[Theorem 8]{MR1839474}
There exists an absolute constant $c>0$ such that for all polynomial $Q:\mathbb{R}^n\rightarrow\mathbb{R}$ of degree at most $d$, all $1<q<\infty$, all log-concave probability measures $\mu$ on $\mathbb{R}^n$  and all $\alpha>0$,
$$\left(\int_{\mathbb{R}^n}|Q(\textbf{x})|^{\frac{q}{d}}\mu(d\textbf{x})\right)^{\frac{1}{q}}\mu(\{\textbf{x}\in\mathbb{R}^n\quad |Q(\textbf{x})|<\alpha\})\leq cq\alpha^{\frac{1}{d}}.$$
\end{lem}
\begin{lem}\label{event} 
Let $1<\beta<2$. There is an event $\Omega^{**}$ of probability one such that for all $\omega\in\Omega^{**}$, 
for all $\epsilon$ in $E$, for all $\textbf{t}_0\in[0,1]^d$, and for $j$ large enough, 
there is a $\textbf{k}_j\in  \mathbb{Z}^d$ such that:
$$|\textbf{t}_0-2^{-j}\textbf{k}_j|\leq bj^{\beta}2^{-j}$$
where $b>0$ is a constant independent of $\omega$ and
$$|I_d(\psi^{\epsilon}_{j,\textbf{k}_j})|>\eta$$
where $\eta>0$ is a constant independent of $\omega$.
\end{lem}
\textbf{Proof}: Let $j$ be a positive integer. For simplicity, we assume $\epsilon=(1,...,1)$.
The general case follows in a similar way. We consider the dyadic points of order $j$ of $(0,1)^d$ {\it i.e.}
$\{\textbf{k}2^{-j}:\textbf{k}\in\{1,...,2^j-1\}^d\}$. Set $\Delta_j=\{1,...,2^j-1\}^d$. Let us build a 
covering of $\Delta_j$. In that view, define the following sets:
$$Q_j=\{0,...,\lfloor \dfrac{2^j}{j^{\beta}}\rfloor+1\}^d,$$
$$\forall \textbf{q}\in Q_j, \ D^{\textbf{q}}_j=\{\lceil j^{\beta}\rceil\textbf{q}+\textbf{r}:\textbf{r}\in\{0,...,\lceil j^{\beta}\rceil-1\}^d\}.$$
Clearly:
$$\Delta_j\subset \bigcup_{\textbf{q}\in Q_j}D^{\textbf{q}}_j.$$
Consider the following event:
$$A_j=\bigcup_{\textbf{q}\in Q_j}\bigcap_{\textbf{k}\in diag(D^q_j)}\{|I_d(\psi^{\epsilon}_{j,\textbf{k}})|<\dfrac{1}{4.d^dc^d}\}.$$
Where $diag(D^q_j)$ denotes the diagonal of $D^q_j$. We wish to estimate its probability, and in particular its dependence on $j$. For this purpose, 
we first perform a study of the random variables $I_d(\psi^{\epsilon}_{j,\textbf{k}})$. 
Since $\epsilon=(1,...,1)$, one has:
$$\psi^{\epsilon}_{j,\textbf{k}}=\psi_{j,k_1}\otimes...\otimes\psi_{j,k_d}$$
and
$$\hat{\psi}^{\epsilon}_{j,\textbf{k}}=\psi_{j,k_1}\hat{\otimes}...\hat{\otimes}\psi_{j,k_d}.$$
Let us reorder the previous symmetric tensor product in the following way:
$$\hat{\psi}^{\epsilon}_{j,\textbf{k}}=\psi_{j,k_1}^{\hat{\otimes}\gamma_1}\hat{\otimes}...\hat{\otimes}\psi_{j,k_l}^{\hat{\otimes}\gamma_l},$$
where $l\in{1,...,d}$, $\gamma_1+...+\gamma_l=d$ and, for all $p,m$ in $\{1,...,l\}$, $k_p\ne k_m$. 
Using the formula expressing Wiener-Itô integrals in terms of Hermite polynomials, one obtains:
$$I_d(\psi^{\epsilon}_{j,\textbf{k}})=\prod_{p=1}^{l}H_{\gamma_p}(\int_{\mathbb{R}}\psi_{j,k_p}(u)dB_u).$$
The standard normal random variables $\int_{\mathbb{R}}\psi_{j,k_p}(u)dB_u$ are independent. 
As a consequence, each $I_d(\psi^{\epsilon}_{j,\textbf{k}})$ is a polynomial of degree at most $d$ 
evaluated at independent standard normal random variables. Moreover, for each $k,k'\in diag(D^q_j)$, 
$I_d(\psi^{\epsilon}_{j,\textbf{k}})$ and $I_d(\psi^{\epsilon}_{j,\textbf{k}'})$ are mutually independent
for $k \neq k'$. Therefore:
$$\mathbb{P}(A_j)\leq \sum_{q\in Q_j}\prod_{\textbf{k}\in diag(D^q_j)}\mathbb{P}(\{|I_d(\psi^{\epsilon}_{j,\textbf{k}})|<\dfrac{1}{4.d^dc^d}\})$$
Lemma \ref{Carb} with $q=2d$ yields
$$\mathbb{P}(\{|I_d(\psi^{\epsilon}_{j,\textbf{k}})|<\dfrac{1}{4.d^dc^d}\})\leq (\dfrac{1}{4})^{\frac{1}{d}}(\sqrt{d!}||\hat{\psi}^{\epsilon}_{j,\textbf{k}}||_2)^{-\frac{1}{d}}.$$
Using the trivial fact:
$$\sqrt{d!}||\hat{\psi}^{\epsilon}_{j,\textbf{k}}||_2\geq 1,$$
one obtains the estimate
$$\mathbb{P}(\{|I_d(\psi^{\epsilon}_{j,\textbf{k}})|<\dfrac{1}{4.d^dc^d}\})\leq (\dfrac{1}{4})^{\frac{1}{d}}$$
which entails that
$$\mathbb{P}(A_j)\leq (\lfloor \dfrac{2^j}{j^{\beta}}\rfloor+2)^d (\dfrac{1}{4})^{\frac{\lceil j^{\beta}\rceil}{d}}.$$
Since $\beta>1$, $\sum_{j=1}^{+\infty}\mathbb{P}(A_j)<\infty$. By the Borel-Cantelli lemma, the following 
event has probability one:
$$\Omega^{**}=\bigcup_{J= 1}^{\infty}\bigcap_{j=J}^{\infty}A^c_j,$$
where $A^c_j$ denotes the complement of the event $A_j$. Let $\omega\in\Omega^{**}$ and 
$\textbf{t}_0\in[0,1]^d$. There exists $\tilde{\textbf{k}}_j\in\Delta_j$ such that $|
\textbf{t}_0-\tilde{\textbf{k}}_j2^{-j}|\leq \sqrt{d}2^{-j}$. 
Moreover, since $\bigcup_{q\in Q_j}D^q_j$ is a covering of $\Delta_j$ there exists $\tilde{q}_j$ such that
$\tilde{\textbf{k}}_j\in D^{\tilde{q}_j}_j$.
But for $j$ large enough, there also exists $\hat{\textbf{k}}_j\in D^{\tilde{q}_j}_j$ such that 
$|I_d(\psi_{j,\hat{\textbf{k}}_j})|>\dfrac{1}{4.d^dc^d}$. Write:
$$|\textbf{t}_0-2^{-j}\hat{\textbf{k}}_j|\leq |\textbf{t}_0-\tilde{\textbf{k}}_j2^{-j}|+2^{-j}|\hat{\textbf{k}}_j-\tilde{\textbf{k}}_j|\leq bj^{\beta}2^{-j},$$
which yields the desired conclusion. $\Box$\\
\\
We are now ready to prove that the pointwise Hölder exponent of the stohastic field $\textbf{X}^{\alpha}$ 
is equal to $H$ almost surely:

\begin{theo}\label{expchamp}
Almost surely,
$$\forall \textbf{t}\in(0,1)^d\quad\gamma_{\textbf{X}^{\alpha}}(\textbf{t})=H.$$
\end{theo}
\textbf{Proof}: Assume that there exists $\omega_0\in\Omega^{**}$, $\textbf{t}_0\in(0,1)^d$ and $\epsilon_0>0$ such that:
$$\gamma_{\textbf{X}^{\alpha}}(\textbf{t}_0)(\omega_0)\geq H+\epsilon_0.$$
By definition of the pointwise Hölder exponent, for any $\textbf{t}$ such that 
$||\textbf{t}-\textbf{t}_0||_2$ is sufficiently small, the following inequality holds:
\begin{equation}\label{abs}
|\textbf{X}^{\alpha}_{\textbf{t}}(\omega_0)-\textbf{X}^{\alpha}_{\textbf{t}_0}(\omega_0)|\leq C(\omega_0)||\textbf{t}-\textbf{t}_0||_2^{H+\epsilon_0},
\end{equation}
where $C(\omega_0)>0$ is a finite constant. By continuity and growth at infinity of 
$\textbf{X}^{\alpha}$, this inequality holds for any $\textbf{t}\in \mathbb{R}^d$. 
For any $(j,\textbf{k},\epsilon)\in\mathbb{N}\times\mathbb{Z}^d\times E$,
$$|I_d(\psi^{\epsilon}_{j,\textbf{k}})(\omega_0)|=2^{j(H+d)}|\int_{\mathbb{R}^d}\textbf{X}^{\alpha}_{\textbf{t}}D^{\alpha+1}(\psi^{\epsilon})(2^j\textbf{t}-\textbf{k})d\textbf{t}|.$$
Using the fact that $D^{\alpha+1}(\psi^{\epsilon})\in S_0(\mathbb{R}^d)$ and inequality \eqref{abs}, 
one gets:
\begin{eqnarray*}
|I_d(\psi^{\epsilon}_{j,\textbf{k}})(\omega_0)| & \leq & 
2^{j(H+d)}\int_{\mathbb{R}^d}|\textbf{X}^{\alpha}_{\textbf{t}}-\textbf{X}^{\alpha}_{\textbf{t}_0}||D^{\alpha+1}(\psi^{\epsilon})(2^j\textbf{t}-k)|d\textbf{t} \\
& \leq & C(\omega_0)2^{j(H+d)}\int_{\mathbb{R}^d}||\textbf{t}-\textbf{t}_0||_2^{H+\epsilon_0}|D^{\alpha+1}(\psi^{\epsilon})(2^j\textbf{t}-\textbf{k})|d\textbf{t}.
\end{eqnarray*}

\noindent The rapid decrease of $D^{\alpha+1}(\psi^{\epsilon})$ and the the change of 
variable $\textbf{x}=2^j\textbf{t}-\textbf{k}$ yield
$$|I_d(\psi^{\epsilon}_{j,\textbf{k}})(\omega_0)|\leq C'(\omega_0)2^{jH}\int_{\mathbb{R}^d}||2^{-j}\textbf{x}+\textbf{k}2^{-j}-\textbf{t}_0||_2^{H+\epsilon_0}\dfrac{1}{(1+||\textbf{x}||_2)^p}d\textbf{x}.$$
Since $(|t|+|s|)^{H+\epsilon_0}\leq C''(|t|^{H+\epsilon_0}+|s|^{H+\epsilon_0})$, one obtains
$$|I_d(\psi^{\epsilon}_{j,\textbf{k}})(\omega_0)|\leq C'''(\omega_0)2^{-j\epsilon_0}(1+||\textbf{k}-2^j\textbf{t}_0||_2^{H+\epsilon_0}).$$
This inequality is true for any $j$ and $\textbf{k}$. Applying it with the $\textbf{k}_j$ of 
Lemma \ref{event} yields
$$|I_d(\psi^{\epsilon}_{j,\textbf{k}})(\omega_0)|\leq C'''(\omega_0)2^{-j\epsilon_0}(1+(bj^{\beta})^{H+\epsilon_0}).$$
This right hand side in the above inequality tends to $0$ as $j$ tends to $+\infty$. 
But, by Lemma \ref{event}, $|I_d(\psi^{\epsilon}_{j,\textbf{k}_j})(\omega_0)|>\eta$. Then, ad absurdum, 
$$a.s, \ \forall \textbf{t}\in(0,1)^d, \ \gamma_{\textbf{X}^{\alpha}}(\textbf{t})=H.\quad \Box$$
\begin{cor}\label{ph} 
Almost surely, 
$$ \forall t\in(0,1), \ \gamma_{X^{\alpha}}(t)=H.$$
\end{cor}
This result should be compared with Proposition 2.2 p.163 and Remark 2.1 p.164 of \cite{Takashima}, which entail that, for all $t$,
almost surely, the pointwise Hölder exponent of $X^{\alpha}$ (as well as the ones of all Hermite processes) are equal to $H$.

\section{Hausdorff Dimensions of SSSI processes represented by multiple Wiener-Itô
integrals}\label{sec:dim}

In this section, we provide estimates for the Hausdorff dimension of the range 
and the graph of multidimensional processes
represented by multiple Wiener-Itô integrals of order $d$ whose kernel 
verify conditions $1.,2.,3.$ with self-similar exponent $H$. 
More precisely, let $d\geq 1$ and $N\geq 1$ be two integers. Let $\{Y^{H}_t\}$ be the process defined by:
$$Y^{H}_t=\gamma(H,d)I_d(h_t^{H})$$
where $\gamma(H,d)$ is a normalizing positive 
constant such that $\mathbb{E}[|Y_1^{H}|^2]=1$ and $h^H_t$ any symmetric square-integrable kernel satisfying properties 1.,2.,3.\\
Let $\frac{1}{2}<H_1\leq...\leq H_{N}<1$. Let 
$\{\mathbb{Y}^{H}_t\}$ be the multidimensional process defined by:
$$\{\mathbb{Y}^{H}_t\}=\{(Y^{H_1}_t,...,Y^{H_{N}}_t): t\in\mathbb{R}_+\}$$
where the coordinates are independent copies of the process $Y^{H}_t$. 

Let $E\subset\mathbb{R}_+$.
It turns out that the methods developed in \cite{MR1415699} apply with few modifications
for estimating the Hausdorff dimension  of the range 
$R_{E}(\mathbb{Y}^{H}):= \{\mathbb{Y}^{H}(t), \ t \in E \}$
and the graph $Gr_{E}(\mathbb{Y}^{H}) := \{(t,\mathbb{Y}^{H}(t)), \ t \in E \}$
of $\mathbb{Y}^{H}$ over $E$. We shall make use of the following lemma which is similar to lemma 4.3 of \cite{NourdinPoly}:

\begin{lem}\label{smb}
There exists a stricly positive constant, $C$, depending only on $d$, such that:
$$\forall j\in\{1,...,N\}\quad \mathbb{P}(|Y^{H_j}_1|\leq x)\leq Cx^{\frac{1}{d}}.$$
\end{lem}
\textbf{Proof}: We shall use Lemma \ref{Carb} in order to prove this inequality. For all 
$j\in \{1,...,N\}$, $\gamma(H_j,d)h_1^{H_j}\in L^2(\mathbb{R}^d)$.
If $\{f_n:n\in\mathbb{N}\}$ is an orthonormal basis of $L^2(\mathbb{R})$, then $\{f_{i_1}\otimes...\otimes f_{i_{d}}:(i_1,...,i_{d})\in\mathbb{N}^{d}\}$ is an orthonormal basis of $L^2(\mathbb{R}^d)$. Thus, one can approximate $\gamma(H_j,d)h_1^{H_j}$ by:
$$h_1^{H_j,n}=\sum_{i_1,...,i_{d}=1}^{n}<\gamma(H_j,d)h_1^{H_j};f_{i_1}\otimes...\otimes f_{i_{d}}>f_{i_1}\otimes...\otimes f_{i_{d}}.$$
By linearity of the multiple Wiener-Itô integral,
$$I_d(h_1^{H_j,n})=\sum_{i_1,...,i_{d}=1}^{n}<\gamma(H_j,d)h_1^{H_j};f_{i_1}\otimes...\otimes f_{i_{d}}>I_d(f_{i_1}\otimes...\otimes f_{i_{d}}).$$
The Wiener-Itô theorem entails that $I_d(f_{i_1}\otimes...\otimes f_{i_{d}})$ is a polynomial of degree 
at most $d$ evaluated at independent standard normal random variables. Therefore,
$$I_d(h^{H_j,n}_1)=P_{n,d}(\int_{\mathbb{R}}f_1(u)dB_j(u),...,\int_{\mathbb{R}}f_n(u)dB_j(u))$$
where $P_{n,d}$ is a polynomial of degree at most $d$. Applying Lemma \ref{Carb} to this quantity 
with $q=2d$ and $\mu$ the Gaussian probability measure on $\mathbb{R}^n$ yields that there exists a 
constant $C(d)>0$ such that:
$$\forall n\in\mathbb{N}\quad\forall x>0\quad\mathbb{P}(|I_d(h^{H_j,n}_1)|\leq x)\leq C(d)x^{\frac{1}{d}}||h^{H_j,n}_1||^{-\frac{1}{d}}_2.$$
By the continuity of the multiple Wiener-Itô integral, $I_d(h^{H_j,n}_1)\rightarrow I_d(\gamma(H_j,d)h^{\alpha_j}_1)$ in $L^2(\Omega)$. Consider a subsequence such that $I_d(h^{H_j,n}_1)$ converges almost surely to $I_d(\gamma(H_j,d)h^{H_j}_1)$. Using Fatou's lemma along this subsequence, one obtains:

\begin{eqnarray*}
\mathbb{P}(I_d(h^{H_j}_1)\leq x) & \leq & \underset{n\rightarrow+\infty}{\liminf} \ \mathbb{P}(|I_d(h^{H_j,n}_1)|\leq x) \\
& \leq & C(d)x^{\frac{1}{d}} \ \underset{n\rightarrow+\infty} {\liminf} \ ||h^{H_j,n}_1||^{-\frac{1}{d}}_2 \\
& \leq & C(d)x^{\frac{1}{d}}\ ||\gamma(H_j,d)h^{H_j}_1||^{-\frac{1}{d}}_2\\
& \leq & C(d)x^{\frac{1}{d}}(\dfrac{1}{\sqrt{d!}})^{-\frac{1}{d}}\\
& \leq & K(d)x^{\frac{1}{d}}
\end{eqnarray*}

since, we have $\mathbb{E}[|Y^{H_j}_1|^2]=\gamma(H_j,d)^2d!||h^{H_j}_1||^2_2=1$.
$\Box$\\
\\
We are now ready to prove the following theorem which is similar to Theorem 3.3 of \cite{MR1415699}.
\begin{theo}\label{lbdim}
Almost surely, 
$$dim_{\mathcal{H}}R_{E}(\mathbb{Y}^{H})\geq \min\left(N;\dfrac{dim_{\mathcal{H}}E+\frac{\sum_{j=1}^k(H_k-H_j)}{d}}{H_k},k=1,...,N\right),$$
$$dim_{\mathcal{H}}Gr_{E}(\mathbb{Y}^{H})\geq \min\left(\dfrac{dim_{\mathcal{H}}E+\frac{\sum_{j=1}^k(H_k-H_j)}{d}}{H_k},k=1,...,N,dim_{\mathcal{H}}E+\sum_{i=1}^{N}\frac{(1-H_i)}{d}\right).$$
\end{theo}
\textbf{Proof}: Assume that $dim_{\mathcal{H}}E >0$ (otherwise there is nothing to prove) 
and let $k$ be such that:
$$\dfrac{\sum_{i=1}^{k-1}H_i}{d}<dim_{\mathcal{H}}E\leq \dfrac{\sum_{i=1}^kH_i}{d}$$
Let us compute the quantity:
$$\lambda=\min\left(N;\dfrac{dim_{\mathcal{H}}E+\frac{\sum_{i=1}^j(H_j-H_i)}{d}}{H_j},j=1,...,N\right).$$
Let $j<k$. Let:
$$\Delta:=\dfrac{dim_{\mathcal{H}}E+\frac{\sum_{i=1}^k(H_k-H_i)}{d}}{H_k}-\dfrac{dim_{\mathcal{H}}E+\frac{\sum_{i=1}^j(H_j-H_i)}{d}}{H_j}.$$
Then,

\begin{eqnarray*}
\Delta & = & \dfrac{1}{H_jH_k}\left[dim_{\mathcal{H}}E(H_j-H_k)+\frac{\sum_{i=1}^k(H_k-H_i)H_j}{d}-\frac{\sum_{i=1}^j(H_j-H_i)H_k}{d}\right] \\
& = & \dfrac{1}{H_jH_k}\left[dim_{\mathcal{H}}E(H_j-H_k)+\frac{\sum_{i=1}^j(H_k-H_j)H_i}{d}+\frac{\sum_{i=j+1}^k(H_k-H_i)H_j}{d}\right].
\end{eqnarray*}

Since:
$$\sum_{i=j+1}^k(H_k-H_i)H_j=\sum_{i=j+1}^{k-1}(H_k-H_i)H_j \leq \sum_{i=j+1}^{k-1}(H_k-H_j)H_i,$$

one gets:
$$\dfrac{dim_{\mathcal{H}}E+\frac{\sum_{i=1}^k(H_k-H_i)}{d}}{H_k}-\dfrac{dim_{\mathcal{H}}E+\frac{\sum_{i=1}^j(H_j-H_i)}{d}}{H_j}\leq\dfrac{H_k-H_j}{H_jH_k}\left[-dim_{\mathcal{H}}E+\frac{\sum_{i=1}^{k-1}H_i}{d}\right]$$
which is strictly less than 0 by definition of $k$. One proves in a similar way that this inequality 
also holds for $j>k$. Thus,
$$\lambda=\dfrac{dim_{\mathcal{H}}E+\frac{\sum_{i=1}^k(H_k-H_i)}{d}}{H_k}.$$
From this one deduces that $\frac{k-1}{d}<\lambda\leq \frac{k}{d}$.\\

We shall use a classical potential theoretic method \cite{MR2118797} to obtain a lower bound on the Hausdorff 
dimension. More precisely, in the case of the range of $\{\mathbb{Y}^{H}_t\}$, 
we study the convergence of the following quantity:
$$\int_{R_{E}(\mathbb{Y}^{H})\times R_{E}(\mathbb{Y}^{H})}\dfrac{d\mu_{\mathbb{Y}^{H}}(x)d\mu_{\mathbb{Y}^{H}}(y)}{|x-y|^{\gamma}}$$
where $\gamma$ is a strictly positive constant and $\mu_{\mathbb{Y}^{H}}$ is the image by $\mathbb{Y}^{H}$ of a positive measure, $\mu$, on $E$, both of them suitably chosen. By a change of variable, the previous quantity is equal to:
$$\int_{E\times E}\dfrac{d\mu(t)d\mu(s)}{|\mathbb{Y}^{H}_t-\mathbb{Y}^{H}_s|^{\gamma}}.$$
By Fubini's theorem, it suffices to show the convergence of the following integral:
$$\int_{E\times E}\mathbb{E}\left[|\mathbb{Y}^{H}_t-\mathbb{Y}^{H}_s|^{-\gamma}\right]d\mu(t)d\mu(s).$$
Fix $\gamma\in (\frac{k-1}{d},\lambda)$ and set
$$F(x)=\prod_{j=1}^{N}\mathbb{P}[|Y^{H_j}_1|\leq x_j].$$
One computes:

\begin{eqnarray*}
\mathbb{E}[|\mathbb{Y}^{H}_t-\mathbb{Y}^{H}_s|^{-\gamma}] & = & \int_{\mathbb{R}_+^{N}}[(x_1|t-s|^{H_1})^2+...+(x_{N}|t-s|^{H_{N}})^2]^{-\frac{\gamma}{2}}dF(x)\\
& = & |t-s|^{-H_1\gamma}\int_{\mathbb{R}_+^{d-1}}dF_2(x_2)...dF_{N}(x_{N})\int_{\mathbb{R}_+}[x_1^2+...+(x_{N}|t-s|^{H_{N}-H_1})^2]^{-\frac{\gamma}{2}}dF_1(x_1).
\end{eqnarray*}

For all $a>0$ and all $\beta>0$:
$$J:=\int_{\mathbb{R}_+}[x_1^2+a^2]^{-\frac{\beta}{2}}dF_1(x_1)=\beta\int_{\mathbb{R}_+}x_1F_1(x_1)[x_1^2+a^2]^{-\frac{\beta}{2}-1}dx_1.$$
When $\beta>\frac{1}{d}$, Lemma \ref{smb} yields:
\begin{eqnarray*}
J & \leq & \beta C(d)\int_{0}^{+\infty}[x_1^2+a^2]^{-\frac{\beta}{2}-1}x_1^{1+\frac{1}{d}}dx_1 \\
& \leq & C(d)\left[a^{-\beta-2}\int_{0}^{a}x_1^{1+\frac{1}{d}}dx_1+\int_{a}^{+\infty}x_1^{-\beta-1+\frac{1}{d}}dx_1\right] \\
& \leq & C'(\beta,d)a^{-\beta+\frac{1}{d}}.
\end{eqnarray*}

\noindent When $\beta<\frac{1}{d}$,
$$\int_{\mathbb{R}_+}x_1F_1(x_1)[x_1^2+a^2]^{-\frac{\beta}{2}-1}dx_1\leq \int_{0}^1x_1^{-\beta-1+\frac{1}{d}}dx_1+\int_{1}^{+\infty}x_1^{-\beta-1}dx_1$$
and thus
$$\int_{\mathbb{R}_+}[x_1^2+a^2]^{-\frac{\beta}{2}}dF_1(x_1)\leq C''(\beta,d).$$
Since $\gamma\in (\frac{k-1}{d},\lambda)$ we see that, upon integrating iteratively with respect to $dF_1(x_1),dF_2(x_2),...,dF_k(x_k)$, the integral
$$\int_{\mathbb{R}_+^{N}}[(x_1|t-s|^{H_1})^2+...+(x_{N}|t-s|^{H_{N}})^2]^{-\frac{\gamma}{2}}dF(x)$$
is not larger than a constant times:
$$|t-s|^{-H_1\gamma-(H_2-H_1)(\gamma-\frac{1}{d})-...-(H_k-H_{k-1})(\gamma-\frac{k-1}{d})}=|t-s|^{-H_k\gamma+\sum_{j=1}^k\frac{H_k-H_j}{d}}.$$
By Frostman's Lemma, the condition $H_k\gamma-\sum_{j=1}^k\frac{H_k-H_j}{d}<dim_{\mathcal{H}}E$ 
entails that there exists a positive measure $\mu$ on $E$  such that
$$\int_{E\times E}\dfrac{\mu(ds)\mu(dt)}{|t-s|^{H_k\gamma-\sum_{j=1}^k\frac{H_k-H_j}{d}}}<+\infty.$$
Therefore,
$$\mathbb{E}\left[\int_{E\times E}\dfrac{\mu(ds)\mu(dt)}{|\mathbb{Y}^{H}_t-\mathbb{Y}^{H}_s|^{\gamma}}\right]<+\infty$$
{\it i.e.}
$$dim_{\mathcal{H}}R_{E}(\mathbb{Y}^{H})\geq \gamma\quad a.s.$$
for all $\gamma\in (\frac{k-1}{d},\lambda)$.\\
\\
The proof of the lower bound for the Hausdorff dimension of the graph of $\mathbb{Y}^{H}_t$ over $E$
follows from the one of Theorem 3.3 in \cite{MR1415699} with the same modifications as above. We omit the
details.\quad $\Box$\\
\\
\textbf{Remarks}:
\begin{itemize}
\item The proof relies essentially on Lemma \ref{smb}. The lower bound holds for multidimensional 
Hermite processes of any order $d$ as well as for multidimensional versions of $X^{\alpha}$.
\item Moreover, it is well-known that the density of the Rosenblatt distribution is continuous and bounded (\cite{MR1650532}, \cite{MR1050053}). As a consequence, $\mathbb{P}(|Y^{H_j}_1|\leq x)\leq Kx$, which yields
the following lower bounds:
$$dim_{\mathcal{H}}R_{E}(\mathbb{Y}^{H})\geq \min\left(N;\dfrac{dim_{\mathcal{H}}E+\sum_{j=1}^k(H_k-H_j)}{H_k},k=1,...,N\right)\quad a.s.$$
$$dim_{\mathcal{H}}Gr_{E}(\mathbb{Y}^{H})\geq \min\left(\dfrac{dim_{\mathcal{H}}E+\sum_{j=1}^k(H_k-H_j)}{H_k},k=1,...,N,dim_{\mathcal{H}}E+\sum_{i=1}^{N}(1-H_i)\right)\quad a.s.$$
One thus retrieves here the result of Theorem 2.3 of \cite{MR1415699}.
\item Actually, for $d=2$, we have the following equality in law:
\begin{align*}
Y^{H_j}_1\overset{(d)}{=}\sum_{n=1}^{\infty}\lambda_{n,H_j}(\eta_n^2-1)
\end{align*}
where $\{\lambda_{n,H_j};n\geq 1\}$ satisfying $|\lambda_{1,H_j}|\geq...\geq|\lambda_{n-1,H_j}|\geq|\lambda_{n,H_j}|\geq...\geq 0$ are the eigenvalues of the Hilbert-Schmidt operator associated with the $L^2(\mathbb{R}^2)$-kernel $h_1^{H_j}$ and $\{\eta_n;n\geq 1\}$ is a collection of independent standard normal random variables. A sufficient condition for $Y^{H_j}_1$ to have a continuous and bounded density is obtained in \cite{Hu} using Malliavin calculus. The boundedness and the continuity of the density rely on the existence of negative moments for $\|D(Y^{H_j}_1)\|_{L^2(\mathbb{R})}$ where $D$ is the Malliavin derivative (see theorem 3.1 of \cite{Hu} and chapter 1.2 of \cite{MR2200233} for a definition of $D$). Lemma 7.1. in \cite{Hu} provides sufficient and necessary conditions for the existence of negative moments for $\|D(Y^{H_j}_1)\|_{L^2(\mathbb{R})}$. It depends on the number of non-zero eigenvalues in $\{\lambda_{n,H_j};n\geq 1\}$. Thus, for such $Y^{H_j}_1$, we have the following lower bounds:
$$dim_{\mathcal{H}}R_{E}(\mathbb{Y}^{H})\geq \min\left(N;\dfrac{dim_{\mathcal{H}}E+\sum_{j=1}^k(H_k-H_j)}{H_k},k=1,...,N\right)\quad a.s.$$
$$dim_{\mathcal{H}}Gr_{E}(\mathbb{Y}^{H})\geq \min\left(\dfrac{dim_{\mathcal{H}}E+\sum_{j=1}^k(H_k-H_j)}{H_k},k=1,...,N,dim_{\mathcal{H}}E+\sum_{i=1}^{N}(1-H_i)\right)\quad a.s.$$
\end{itemize}
Upper bounds for the Hausdorff dimensions of the range and the graph are obtained by a direct application
of Theorem 3.1 in \cite{MR1415699}, since we deal with elements of Wiener's chaoses, 
which thus have finite moments of all order:
\begin{theo}\label{ubdim}\cite[Thoerem 3.1]{MR1415699}
Almost surely,
$$dim_{\mathcal{H}}R_{E}(\mathbb{Y}^{H})\leq \min\left(N;\dfrac{dim_{\mathcal{H}}E+\sum_{j=1}^k(H_k-H_j)}{H_k},k=1,...,N\right),$$
$$dim_{\mathcal{H}}Gr_{E}(\mathbb{Y}^{H})\leq \min\left(\dfrac{dim_{\mathcal{H}}E+\sum_{j=1}^k(H_k-H_j)}{H_k},k=1,...,N ; dim_{\mathcal{H}}E+\sum_{i=1}^{N}(1-H_i)\right).$$
\end{theo}
\textbf{Remarks}:
\begin{itemize}
\item For the $N$-dimensional anisotropic Rosenblatt process, we retrieve the known results:
$$dim_{\mathcal{H}}R_{E}(\mathbb{Y}^{H})= \min\left(N;\dfrac{dim_{\mathcal{H}}E+\sum_{j=1}^k(H_k-H_j)}{H_k},k=1,...,N\right)\quad a.s.,$$
$$dim_{\mathcal{H}}Gr_{E}(\mathbb{Y}^{H})= \min\left(\dfrac{dim_{\mathcal{H}}E+\sum_{j=1}^k(H_k-H_j)}{H_k},k=1,...,N ; dim_{\mathcal{H}}E+\sum_{i=1}^{N}(1-H_i)\right)\quad a.s.$$
\item In the isotropic case, the first equality was proven in \cite{MR2759163}.
\end{itemize}

\textbf{Acknowledgements}\\
\\
The author would like to thank his supervisor Jacques Lévy-Véhel for his useful comments and for stimulating discussions which have
helped to greatly improve this paper.

\def\cprime{$'$}

\end{document}